\documentclass[10pt]{amsart}
\usepackage{latexsym}
\usepackage{amssymb}
\usepackage{amsfonts}
\usepackage{epsfig}
\usepackage{amsmath}
\usepackage{amscd}
\usepackage{psfrag}
\usepackage[all]{xy}
\newtheorem{defn0}{Definition}[section]
\newtheorem{prop0}[defn0]{Proposition}
\newtheorem{thm0}[defn0]{Theorem}
\newtheorem{lemma0}[defn0]{Lemma}
\newtheorem{corollary0}[defn0]{Corollary}
\newtheorem{example0}[defn0]{Example}
\newtheorem{conjecture0}[defn0]{Conjecture}
\newtheorem{notation0}[defn0]{Notation}

\theoremstyle{remark}
\newtheorem{remark0}[defn0]{Remark}

\newenvironment{prop}{\begin{prop0}}{\end{prop0}}
\newenvironment{thm}{\begin{thm0}}{\end{thm0}}
\newenvironment{lem}{\begin{lemma0}}{\end{lemma0}}
\newenvironment{cor}{\begin{corollary0}}{\end{corollary0}}
\newenvironment{example}{\begin{example0}\rm}{\end{example0}}
\newenvironment{rem}{\begin{remark0}\rm}{\end{remark0}}

\renewenvironment{proof}{\noindent {\textsc{Proof.}}}{$\square$ \vspace{3mm}}
\newenvironment{proofof}[1]{\noindent {\textsc{Proof of {#1}.}}}{$\square$ \vspace{3mm}}

\newcommand{\I}{\mathbf{I}_{R}}
\newcommand{\LL}{{\mathcal L}}
\newcommand{\K}{\mathcal K}
\newcommand{\G}{\Gamma}

\newcommand{\T}{{\mathcal T}}

\newcommand{\B}{{\mathcal B}}
\newcommand{\m}{\mathfrak{m}}

\newcommand{\OO}{{\mathcal O}}

\newcommand{\br}{\mathfrak{br}_{(X,Q)}}
\renewcommand{\v}{\omega}

\newif\ifprivate
\privatetrue

\def\???{\ifprivate {\bf {???}} \marginpar{{\Huge {\bf ?}}}
\else \fi}

 \DeclareMathOperator{\mult}{mult}
\DeclareMathOperator{\ch}{ch}

  \DeclareMathOperator{\emdim}{emdim}
 \DeclareMathOperator{\Spec}{Spec}

\begin{document}

\title[Exceptional locus on normal surface singularities into a plane]
{On the exceptional locus of the birational projections of a normal surface singularity into a plane}

\author{Jes\'{u}s Fern\'{a}ndez-S\'{a}nchez}
\thanks{This research has been partially supported by the Spanish Committee for Science and Technology (DGYCIT),
projects MTM2005-01518 and MTM2006-14234-C02-02, and the Catalan Research Commission. The author
completed this work as a researcher of the program Juan de la Cierva.
\newline{AMS 2000 subject classification 14B05; 14E05; 14J17}}

\address{Departament de Matem\`atica Aplicada I \\ Universitat Polit\`ecnica de Catalunya \\
Av. Diagonal 647, 08028-Barcelona, Spain.}
 \email{jesus.fernandez.sanchez@upc.edu}

\vspace{6mm} \begin{abstract} 
Given a normal surface singularity $(X,Q)$ and a birational morphism to a non-singular surface $\pi:X\rightarrow S$, 
we investigate the local geometry of the exceptional divisor $L$ of $\pi$. 
We prove that the dimension of the tangent space to $L$ at $Q$ equals the number of exceptional components meeting at $Q$. 
Consequences relative to the existence of such birational projections contracting a prescribed number of irreducible curves are deduced. A new characterization of minimal singularities is obtained in these terms. 
\end{abstract}

\maketitle

\section{Introduction}
Normal surface singularities that can be projected birationally to a non-singular surface are known as \emph{sandwiched} singularities. 
They are rational singularities and among them, are included all
cyclic quotients and minimal surface singularities.
Given the germ of a sandwiched singularity $(X,Q)$, there
exist several birational morphisms from it to a
non-singular surface.
Any such birational projection $\pi:(X,Q)\rightarrow (S,O)$ determines a complete $\m_{S,O}$-primary ideal $I\subset {\mathcal O}_{S,O}$ so that $X$ is the surface obtained by blowing-up $S$ along $I$ (see \cite{Spivak} for details). 
Such a birational projection determines also (and is in turn determined by) an exceptional curve on $(X,Q)$ (see \cite{Le}). 

The main purpose of the present paper is the study of the local geometry of the exceptional curves of such projections. 
We give necessary conditions for a curve on $(X,Q)$ to be the exceptional locus of a birational projection to a non-singular surface. These results complement in some sense those of \cite{AlbFer}, where the topological types of the ideals determined by such birational projections were characterized. 
We will make use of the theory of clusters of infinitely near points and its connection with the theory of complete ideals in a regular local ring of dimension two (see \cite{CasasBook,Lip94a}) and the study of sandwiched singularities (see \cite{moi1,AlbFer}). 

The organisation of the paper is as follows. 
In section 2, we recall definitions and prove technical facts about infinitely near points and sandwiched surface singularities. 
Let $(R,\m)$ be a regular local ring of dimension two over $\mathbb{C}$ and write $S=\Spec(R)$. If $I$ is a complete $\m$-primary ideal in $R$, write $X=Bl_{I}(S)$ for the surface obtained by blowing-up $I$. 
In section 3 and given a singularity $Q\in X$, we prove the existence of principal curves on $X$ with prescribed intersection multiplicities with the exceptional components on $(X,Q)$ (Theorem \ref{simpint}) and we provide an effective procedure to explicitly determine these curves. 
Then we derive consequences. First, we show that the dimension of the tangent space to the exceptional divisor of $X$ at $Q$ is maximal (Theorem \ref{tangent_space}) i.e. it equals the number of its irreducible components  meeting at $Q$.
It also follows that the reduced exceptional curve at $Q$ has only smooth branches with independent tangencies, so it has only minimal singularities (Corollary \ref{excep_minimal}).
In particular, we have that for any birational projection $\pi$ of $(X,Q)$ to a non-singular surface, the number $\mathfrak{l}_{(X,Q)}(\pi)$ of branches at $(X,Q)$ being contracted by $\pi$ satisfies 
\begin{eqnarray}
 \label{intro_aux} \mathfrak{l}_{(X,Q)}(\pi)\leq \emdim_Q(X)=\mult_Q(X)+1.
\end{eqnarray}
In section 4 we give a formula for the number $\br$ of branches of a generic hypersurface section through $(X,Q)$ in terms of the base points of $I$ (Proposition \ref{branches}). 
From it, we show that there exists no birational projections of $(X,Q)$ to a non-singular surface contracting more than $\br+1$ branches through $(X,Q)$ (Corollary \ref{bound_nQ}). 
This fact implies a second bound
\begin{eqnarray*}
 \label{intro_aux2}\mathfrak{l}_{(X,Q)}(\pi) \leq \br+1
\end{eqnarray*}
which is sharper than (\ref{intro_aux}).
%
Finally, we prove that minimal singularities are those normal surface singularities that admit a birational projection $\pi$ for which the bound (\ref{intro_aux}) is attained (Theorem \ref{nBP1}). 

\section{Preliminaries}
Throughout this work the base field is the field $\mathbb{C}$ of complex
numbers. A \emph{curve} will always be an effective Weil divisor on a surface. We use the symbol $\sharp$ as meaning \emph{cardinality}.

\subsection{Clusters of base points of complete ideals}
A reference for some of the material treated here
is the book \cite{CasasBook} and the reader
is referred to it for proofs. Let $S$ be a regular surface
over $\mathbb{C}$ and $O\in S$. Write $(R,\m)$
for the local ring ${\mathcal O}_{S,O}$. A
\emph{cluster} of points with origin $O$ is a finite set
$K$ of points infinitely near or equal to $O$ such that for
any $p\in K$, $K$ contains all points preceding $p$. 
A system of virtual multiplicities for a cluster $K$ is a collection of integers $\nu=\{\nu_p\}$. The pair ${\mathcal
K}=(K,\nu)$ is called a \emph{weighted cluster}. We write $p\geq q$
if $p$ is infinitely near or equal to $q$, and
$p\rightarrow q$ if $p$ is proximate to $q$. If $p$ is
maximal among the points of $K$ proximate to $q$, then we
say that $p$ is $m_K$-\emph{proximate} to $q$, and denote
it by $p\rightarrow_{m_K}q$; $p$ is said to be
$m_K$-\emph{free} or $m_K$-\emph{satellite} according to
if it is $m_K$-proximate to one or two points. The
\emph{excess} of $\K$ at $p$ is $\rho^{\K}_p=\nu_p-
\sum_{q \rightarrow p} {\nu _q}$,
and \emph{consistent} clusters are those clusters with no negative excesses;  $\K_+=\{p\in K|\rho^{\K}_p>0\}$ is the set of \emph{dicritical}
points of $\K$. 

We will denote by $\I$ the semigroup of complete $\m$-primary deals in $R$. If $\K$ is a weighted cluster, the equations of all curves going
through it define an ideal $H_\K\in \I$ (see
\cite{CasasBook} \S 8.3).
Any ideal $J\in \I$ has a cluster of base points, denoted
by  $BP(J)$, which consists of the points shared by, and the multiplicities of,
the curves defined by a generic element of $J$. Moreover, the maps $J\mapsto
BP(J)$ and $\K\mapsto H_\K$ are reciprocal isomorphisms
between $\I$ and the semigroup of consistent clusters with no points with virtual multiplicity zero (see \cite{CasasBook} 8.4.11 for details). If $p$ is
infinitely near or equal to $O$, $I_p$ is the simple (i.e. irreducible complete) ideal generated by the
equations of the branches going through $p$, and $\K(p)$ is the
weighted cluster corresponding to it by the above isomorphism. 
Moreover, $\{I_p\}_{p\in \K_+}$ is just the set of the simple ideals appearing in the factorisation of $I$; indeed, $I=\prod_{p\in \K_+}I_p^{\rho_p^{\K}}$ (Theorem 8.4.8 of \cite{CasasBook}).
Consistent clusters are characterised as those clusters whose virtual multiplicities can
be realized effectively by some curve on $S$. If $\K$ is not
consistent, $\widetilde{\K}$ is the cluster obtained from ${\mathcal
K}$ by {\em unloading}, i.e. $\widetilde\K$ is the unique consistent
cluster having the same points as $K$ and such that $H_{\widetilde{\K}}=H_{\K}$ (cf.
\cite{CasasBook} \S 4.2 and \S 4.6).

If $\pi_K:S_K\longrightarrow S$ is the composition of the
blowing-ups of all points in $K$, write $E_K$ for the
exceptional divisor of $\pi_K$ and $\{E_p\}_{p\in K}$ for
its irreducible components. If $C$ is a curve on $S$, $e_p(C)$ is the multiplicity of $C$ at $p$ and 
$v_p(C)$ is the value of $C$ relative to the divisorial
valuation associated to $E_p$. 
Use
$|\,\,{\cdot}\,\,|$ as meaning intersection number and
$[\,\,,\,\,]_P$ as intersection multiplicity at $P$. We
have the equality (projection formula) for $\pi_K$:
$|\pi_K^*{C}{\cdot}D|_{S_K}=[C,(\pi_K)_*D]_O$, $D$ being a
curve on $S_K$ without exceptional components. If $p\in K$,
$|E_p{\cdot}E_p|_{S_K}=-r_p-1$, where $r_p$
is the number of points in $K$ proximate to $p$. 
If $\widetilde{C}^K$ is the strict transform of $C$ on $S_K$, we have that for all $p\in K$
\begin{eqnarray}\label{for:int/excep}
|\widetilde{C}^K{\cdot}E_p|_{S_K}=e_p(C)-\sum_{q\in K,q
\rightarrow p}e_q(C)
\end{eqnarray}

If $K$ is a  (non-weighted) cluster, $\G_K$ will be its dual graph (\S
4.4 of \cite{CasasBook}). 
Unless some confusion may arise, we will
identify the points with the corresponding vertices in $\G_K$. 
Given two
points $q,p$ in $K$, the chain $\ch(q,p)$ is
the linear subgraph of $\G_K$ consisting of all vertices and edges
between $q$ and $p$; $\ch^0(q,p)=\ch(q,p)\setminus
\{q,p\}$. 
We say that $p,q\in K$ are \emph{adjacent} if the vertices associated to them in $\G_K$ are connected by one edge. 

\begin{lem}[Lemma 1.3 of \cite{AlbFer}]
\label{chain} Let $p,q\in K$ and 
write $\ch(q,p)=\{u_0=q,u_1,\ldots,u_n,u_{n+1}=p\}$. If $p$ is infinitely near to $q$, then there exists
some $i_0\in \{0,\ldots,n+1\}$ such that
\begin{eqnarray*}
\begin{array} {lc} u_k \leftarrow
u_{k+1} & \mbox{ if }k\in \{0,\ldots i_0-1\} \\ u_k \rightarrow u_{k+1} &
\mbox{ if }k\in \{i_0,\ldots,n\}.
\end{array}
\end{eqnarray*}
Furthermore, if $j\geq i_0$, $u_j$ is proximate to some $u_{\sigma(j)}$ with
$\sigma(j)\leq i_0-1$.
\end{lem}

Let $\K=(K,\nu)$ be a consistent cluster. 
Write $\K_{\v}$ for the cluster obtained from
$\K$ by adding some $\v\in E_K$ as a point with virtual multiplicity one. Then, 
the ideal $H_{\K_{\v}}\subset I$ has codimension one in
$I$, and every complete $\m$-primary ideal
of codimension one in $I$ has this form for some ${\v}
\in E_K$ (Lemma 3.1 of \cite{Edinburgh}). If $\K_{\v}$ is
not consistent, write $\K'=(K',\nu')$ for the cluster obtained from $\widetilde{\K_{\v}}$ by dropping the points with virtual multiplicity 0. Write also 
\begin{eqnarray}
\label{defn:Tq} T_{\v}=\{p\in K\mid
v_p^{\K'}>v_p^{\K}\}.
\end{eqnarray}
We have that $T_{\v}\cap \K_+=\emptyset$. All the unloading steps leading from
$\K_{\v}$ to $\widetilde{\K_{\v}}$ are tame (Remark 4.2 of
\cite{moi1}), so each unloading on a point, say $p$, increases by one the virtual value on $p$
while that of the other points remain unaffected. $T_{\v}$ is just the set of points where some unloading is performed.
There exists a unique minimal point in $T_{\v}$, that we will denote $o_{\v}$. For this point,  $\nu'_{o_{\v}}=\nu_{o_{\v}}+1$, while for $p\neq o_{\v}$, we have $\nu_p-1 \leq \nu'_p \leq \nu_p$ (see Lemma 3.9 of \cite{Edinburgh}). Write \[\B_{\v}=\{p\in K \mid \nu'_p=\nu_p-1\}\]and for
each $p\in K$, $\varepsilon_p=\nu'_p-\nu_p$. 
Write also $\K_+^{\v}=\{p\in \K_+ \mid p \mbox{
adjacent to some point of }T_w\}$.

\begin{lem} \label{difexcess}
Write $\rho'_p$ for the excess of $\K'$ at $p$. Then
\begin{itemize}
\item[(a)] $\rho'_p=\rho_p+\varepsilon_p-\sum_{q\rightarrow
p} \varepsilon_q$.
\item[(b)]  If $p\in K$, then $\rho'_p\geq \rho_p$ if $p\in T_{\v}$; $\rho'_p=\rho_p-1$ if $p\in \K_+^{\v}$;   $\rho'_p=\rho_p$, otherwise.
\end{itemize}
\end{lem}

\begin{proof}
(a) is Lemma 3.11 of \cite{Edinburgh}. (b) is a consequence of the proof of Theorem 4.1 of \cite{Edinburgh}.
\end{proof}

\subsection{Sandwiched surface singularities}
The reader is referred to \cite{Spivak,moi1} for proofs and known facts
about sandwiched singularities and their relation with complete ideals. Let $I\in \I$ and write $\K=BP(I)$ for the cluster of its base points
and $\pi:~X=Bl_I(S)\longrightarrow S$ for the blowing-up of $I$. There is
a commutative diagram
\[\xymatrix{  {S_K} \ar[r]^{f}\ar[rd]_{\pi_K} & {X} \ar[d]^{\pi} \\ &  {S}}\]
where the morphism $f$, given by the universal property of
the blowing-up, is the minimal resolution of the
singularities of $X$ (Remark I.1.4 of \cite{Spivak}). These
singularities are by definition sandwiched singularities. 
There is a bijection between the set of simple ideals
$\{I_p\}_{p\in \K_+}$ in the (Zariski) factorisation of $I$
and the set of irreducible components of $\pi^{-1}(O)$
(see Corollary I.1.5 of \cite{Spivak}). We write
$\{L_p\}_{p\in \K_+}$ for the set of these components.
If $C$ is a curve on $S$, we write $\widetilde{C}$ for the strict transform of $C$ on $X$ and $\LL_C=\sum_{p\in \K_+}v_p(C)L_p$. Note that $\widetilde{C}+\LL_C$ is the total transform of $C$ on $X$. If $p\in \K_+$, we will also write $\LL_p$ for $\LL_{C_p}$ where $C_p$ is a curve defined by a generic element of $I_p$. %

For any $Q$ in the exceptional locus of $X$, write ${\mathcal M}_Q$ for the
ideal sheaf of $Q$ in $X$. Then, the ideal
$I_Q:=\pi_*({\mathcal M}_QI{\mathcal O}_X)\subset I$ is
complete, $\m$-primary and has codimension one
in $I$. In fact, the map $Q\mapsto I_Q$ defines a bijection
(Theorem 3.5 of \cite{moi1}):
 \begin{eqnarray}
 \label{bij}
\left\{ \begin{array}{c}
 \mbox{points in the exceptional} \\
 \mbox{locus of $X$} \\
 \end{array} \right\} & \longleftrightarrow &
\left\{ \begin{array}{c}
 \mbox{complete $\m$-primary ideals} \\
 \mbox{of codimension one in $I$} \\
 \end{array} \right\}
\end{eqnarray}
For any $Q$ in the exceptional locus of $X$, there exists some ${\v}\in E_K$ such that $I_Q=H_{\K_{\v}}$. Moreover, $Q$ is singular if and only if the cluster $\K_{\v}$ is not consistent (Proposition
4.4 of \cite{moi1}). 
Once a singular point $Q$ in $X$ has been fixed, a number of objects are attached to it: if $I_Q=H_{\K_{\v}}$ we write $T_Q$, $\K_+^Q$, $o_Q$, $\B_Q$ for $T_{\v}$, $\K_+^{\v}$, $o_{\v}$, $\B_{\v}$, respectively, and all of them are well defined. We will also write $\K_Q$ for the cluster of base points of $I_Q$. 
Note that $\{E_p\}_{p\in T_Q}$ is the set of exceptional components on $S_K$ contracting by $f$ to $Q$ and $\{L_p\}_{p\in \K_+^Q}$ is the set of exceptional components on $X$ meeting at~$Q$.  
%

\subsection{Some technical results}
Let $\v \in E_K$ such that $\K_{\v}$ is not consistent. The aim here is to give some results relating the structure of the graph $\G_K$ and the excesses of the cluster $\K'=\widetilde{\K_{\v}}$. They will be repeatedly used in forthcoming sections. 
For each $p\in K$, write $\overline{z}_p=v_p(H_{\K_{\v}})-v_p(I)$. 
\begin{rem}
Let $Q\in X$ be the singularity associated to $H_{\K_{\v}}$ by the bijection (\ref{bij}). Denote by $Z_Q=\sum_{p\in T_Q}z_p E_p$ the fundamental cycle of $Q$ (see \cite{Art66}). By virtue of Corollary 3.6 of \cite{moi1}, we know that for each $p\in K$, $\overline{z}_p=z_{p}$ if $p\in T_{Q}$ and $0$ otherwise. The following lemma provides a method for computing the coefficients of $Z_{Q}$ from the proximities between the points of $K$ (cf. \cite{Lau71}). However, we state it here independently of its interpretation in terms of sandwiched singularities. 
\end{rem}

\begin{lem}
\label{coef_fund}
We have that 
\begin{itemize}
\item[(a)]
$\overline{z}_p=\varepsilon_p+\sum_{p\rightarrow q}
\overline{z}_q$.
\item[(b)] If $p\in T_{\v}$ verifies one of the following conditions:
\begin{itemize}
\item[(i)] $p=o_{\v}$; \item[(ii)] $\K$ has positive excess at some
point proximate to $p$; \item[(iii)] $p$ is proximate to
some point not in $T_{\v}$;
\end{itemize} then $\overline{z}_p=1$.
\item[(c)] If $u\notin T_{\v}$, then $u$ is proximate to some point in $T_{\v}$ if and only if $u\in \B_{\v}$.
\end{itemize}
\end{lem}

\begin{proof}
In p. 141 of \cite{CasasBook}, it is shown that $v_p(C)=e_p(C)+\sum_{p\rightarrow q}v_q(C)$. Thus, for all $p\in K$, 
$\overline{z}_p = (\nu'_p-\nu_p)+\sum_{p\rightarrow
q}(v'_q-v_q)= \varepsilon_p+\sum_{p\rightarrow
q}\overline{z}_q$. This proves (a).
From this, (b) follows easily. (c) follows from the statement (a) and (\ref{defn:Tq}) using the fact that $q\in T_{Q}$ if and only if $\overline{z}_{q}\geq 1$.
\end{proof}

The following lemma will be needed in the following section.

\begin{lem} \label{key/prop}
\begin{itemize}
\item[(a)]Let $p,q\in K$ with $p$ infinitely near to $q$ and keep the notation
of \ref{chain}. Assume that $\rho'_{u_j}=0$ if $j\in
\{0,1,\ldots,\sigma(n+1)-1,\sigma(n+1)\}$. If $p\in
B_{\omega}$, then $u_j\in \B_{\omega}$ for all $j\in
\{0,1,\ldots,\sigma(n+1)-1\}$. In particular, there exists
some $u\in \ch(q,p)$ with $u\in \B_{\omega}$ and proximate
to $q$.

 \item[(b)]Let $u\in T_{\omega}$ and $p_1,p_2,p_3\in \K_+^{\omega}$ such that $\ch(u,p_i)\cap \ch(u,p_j)=\{u\}$ if $i\neq j$. Assume that
$\rho'_v=0$ for each $v\in \ch^0(u,p_i)$ and $i\in
\{1,2,3\}$. Then, $\rho'_u\geq 2$.
\end{itemize}
\end{lem}

\begin{proof}
From (a) of
\ref{difexcess} and the assumption above, we see that each
$u_j$ above is in $\B_{\omega}$. This proves (a).
To prove (b), notice that any
$w\in \ch(p_i,u)$ with $w\neq p_i$ is in $T_Q$. We distinguish different cases
according to the number of $p_i$'s which are infinitely near to $u$.

\vspace{2mm} \textsc{Case 1} Assume that $p_1,p_2$ and
$p_3$ are infinitely near to $u$. Then, by 
\ref{chain}, each $p_i$ is proximate to some point of $T_Q$
and by (c) of \ref{coef_fund}, $p_i\in \B_{\omega}$. (a)
applies to show that there is some $u^i\in \ch(u,p_i)$ in
$\B_Q$ proximate to $u$. By (a) of \ref{difexcess},
$\rho'_u\geq 3$ if $u\notin \B_{\omega}$ and $\rho'_u\geq
2$, if $u\in \B_{\omega}$. In any case, the claim follows.

\vspace{2mm} \textsc{Case 2} Now we deal with the case
where there is at least some $p_i$ which is not infinitely
near to $u$. In fact, we are showing that  this case cannot
occur with the above assumptions. Write $x_i$ for the
maximal point such that both $p_i$ and $u$ are infinitely
near to it. Note that $x_i\in \ch(p_i,u)$.

If $p_1,p_2$ are infinitely near to $u$, the same argument
used in Case 1 shows that $\rho'_u\geq 1$, the equality
holding if and only if $u\in \B_{\omega}$. In this case,
(a) shows that there are points $w\in \ch(x_3,u)$ and $w'\in
\ch(x_3,p_3)$, both in $\B_{\omega}$ and proximate to $u$.
By (a) of \ref{difexcess}, $\rho'_{x_3}\geq 1$, against the
assumption.

Now, assume that $p_2,p_3$ are not infinitely near to $u$.
In this case, $u$ is infinitely near or equal to $x_2$ and
$x_3$ and so, they are in the same branch of $K$. Since
$\ch(u,p_2)\cap \ch(u,p_3)=\{u\}$, we have that $x_2\neq
x_3$. We can assume that $x_3$ is infinitely near or equal
to $x_2$. By (c) of \ref{coef_fund}, $p_3\in \B_{\omega}$ and since
$\rho'_{x_3}=0$, (a) applies to show that $x_3\in
\B_{\omega}$ and there is some $w\in \ch(x_2,p_2)$ in
$\B_{\omega}$ proximate to $x_2$. If $x_2=p_2$, this leads
to contradiction with (c) of \ref{coef_fund}. If $p_2\neq x_2$, then
$p_2\in \B_{\omega}$ by (c) of \ref{coef_fund}, and by (a) again,
there is some $w'\in \ch(x_2,p_2)$ in $\B_{\omega}$ and
proximate to $x_2$. By applying (a) of \ref{difexcess}, we
see that $\rho'_{x_2}>0$ against the assumption.
\end{proof}

Finally, we state the following proposition.
\begin{prop}
\label{excessinchains} Assume that $\K_{\omega}$ is not consistent. Let
$p,q\in \K_+^{\omega}$. Then, there exists some $u\in
\ch^0(p,q)$ such that $\rho'_u>0$.
\end{prop}

\begin{proof}
If $p$ is infinitely near to $q$, $p\in \B_{\omega}$ by (c) of \ref{coef_fund}
and the claim follows from (a) of
\ref{key/prop} and (c) of \ref{coef_fund}. If $q$ is infinitely
near to $p$, the same argument works. Otherwise, write
$x_0\in K$ for the maximal point such that both $p$ and $q$
are infinitely near to it. Then, $x_0\in \ch(q,p)$. If
$\rho'_{u}=0$ for all $u\in \ch^0(q,p)$, (a) of
\ref{key/prop} applies to $\ch(x_0,q)$ and $\ch(x_0,p)$
and there are $w\in \ch(x_0,q)$, $w'\in \ch(x_0,p)$
such that $w,w'\in \B_{\omega}$ and proximate to $x_0$. By
(a) of \ref{difexcess}, we see that $\rho'_{x_0}>0$ against
the assumption.
\end{proof}

\section{Cartier divisors with prescribed intersection multiplicities}
\begin{thm}\label{simpint}
Let $Q$ be a singular point on $X$ and for each $p\in \K_+^{Q}$, let
$\alpha_p\in \mathbb{Z}_{>0}$. There exists a curve $C$ on $S$ such that the strict transform $\widetilde{C}$ on $X$ is a Cartier divisor that
intersects the exceptional locus of $X$ only at $Q$ and
$[\widetilde{C},L_p]_Q=\alpha_p$, for all $p\in \K_+^{Q}$.
\end{thm}

First, we need the following lemma. 

\begin{lem}\label{Cart}
Let $\{\alpha_p\}_{p\in \K_+}$ be strictly positive integers and let $C\subset S$ be a curve.
If ${\LL_C=\sum_{p\in \K_+}\alpha_p\LL_p}$, then $\widetilde{C}$ is Cartier and $|\widetilde{C}\cdot L_p|_X=\alpha_p$ for all $p\in K$.
\end{lem}

\begin{proof}
Let $g:X'\rightarrow X$ be a resolution of $X$. By applying
the projection formula for $\pi_K$ and for $g$, we obtain that 
$|\widetilde{C}\cdot L_u|_X= -|\LL_C\cdot L_u|_X, \mbox{ for  }u \in \K_+.$
In particular, if we apply this to a generic curve going
through $\K(p)$ ($p\in \K_+$), we have
$|\mathcal{L}_p\cdot L_u|_X=-1$ if $p=u$
and $0$ otherwise. From this, the second claim follows. The first one is a consequence of Theorem 3.1 of \cite{moi4}.
\end{proof}

\begin{proofof}{\ref{simpint}} 
First, we proceed to explain the idea of the proof. Write
$\K_+^{Q}=\{p_1,\ldots,p_s\}$ and $J=\prod_{i=1}^s
I_{p_i}^{\alpha_{p_i}}$. For each $i\in \{1,\ldots,s\}$,
write $s_i$ for the only point in $T_Q$ adjacent to
$p_i$. Write $\K_0$ for the cluster of base points of
$J$. 
In order to prove the existence of curves with the
desired properties, we shall construct a sequence of
(consistent) clusters
\begin{eqnarray*}
\K_0,\K_1, \ldots , \K_n, \quad
\K_j=(K_j,\nu_j)
\end{eqnarray*}
so that 
\textbf{1.} for each $j\in \{1,\ldots,n\}$, $\LL_{\K_j}=\LL_{\K_0}=\sum_{i=1}^s\alpha_{p_i}\LL_{p_i}$.
\textbf{2.} the dicritical points of $\K_n$ are within $T_Q$.
Then, we take $C$ to be a generic curve going through $\K_n$, so that $\LL_C=\sum_{i=1}^s\alpha_{p_i}\LL_{p_i}$. Condition \textbf{2.} ensures that $\widetilde{C}$ meets the exceptional locus of $X$ only at $Q$, while \textbf{1.} ensures the remaining desired properties by virtue of \ref{Cart}.

Let $w(0)$ be any point in $E_Q$, and set $\K_1$ for the cluster
obtained by adding $w(0)$ as a
simple point to $\K_0$, unloading multiplicities and dropping the points with virtual
multiplicity zero. By (b) of \ref{difexcess}, we have that $\rho^{\K_1}_{p_i}=\alpha_{p_i}-1\geq 0$, $i=1,\ldots,s$. 
For $k\geq 1$ and as far as there exists
some $i\in \{1,\ldots,s\}$ such that $\rho^{\K_k}_{p_r}>0$,
define $\K_{k+1}$ as follows:
write \[T_{w(k)}=\{p\in K_k \mid v_p^{\K_{k}}>
v_p^{\K_{k-1}}\}\] and note that
$\K_k$ has excess $0$ at the points of $T_{w(k)}$ (see \S 2.2). Choose
any $p_r$ with $\rho^k_{p_r}>0$ as above and write $w(k+1)$ for the minimal point
not in $K_k$ which is proximate to $p_r$ and infinitely
near to $s_r$; take $\K_{k+1}$ as the consistent
cluster obtained from $\K_k$ by adding $w(k+1)$, unloading
multiplicities and dropping 
the points with virtual multiplicity zero.

I claim that after finitely many steps, we
reach some $\K_n$ such that $\rho^{\K_n}_{p_i}=0$ for each
$i\in \{1,\ldots, s\}$.  To prove this, we must show that each step above does not decrease the excess of any $p_i$ if $i\neq r$. This is a direct consequence of the following lemma. 

\begin{lem}
If $k\geq 1$ and $i\neq j$, there exists some $u\in
\ch^0(p_i,p_j)$ such that $\rho^k_u>0$.
\end{lem}

\begin{proof}
We use induction on $k$. The case $k=1$ is the claim of 
\ref{excessinchains}. Assume that $\rho^k_{p_r}>0$ and $\K_{k}$ is obtained by adding a point $w(k)$ proximate to $p_r$ as
described above. By \ref{excessinchains} and the induction hypothesis, it is clear that for any $i\neq r$, there is some $u\in \ch^0(p_r,p_i)$ such that $\rho^k_u=0$. Now, take $r\notin \{i,j\}$ and assume that there is no $u\in \ch^0(p_i,p_j)$ with $\rho^k_u>0$. Then, by the induction hypothesis, there is some $u'\in \ch^0(p_i,p_j)$ such that $\rho^{k-1}_{u'}>0$ and necessarily, $\rho^k_{u'}=\rho^{k-1}_{u'}-1$, while $\rho^{k-1}_{u}=0$ if $u\neq u'$. In this case, apply (b) of \ref{key/prop} to deduce that $\rho^k_{u'}\geq 2$ and hence, $\rho^k_{u'}>0$ against the assumption. 
\end{proof}

Therefore, after finitely many steps we get a 
cluster $\K_n$ such that $\rho^{\K_n}_{p}=0$ if $p\notin T_Q$. This gives condition \textbf{2}. Since no unloading step is performed throughout the above procedure on any $p\in \K_+$, condition \textbf{1} above is also satisfied. 
Define $\T=(K',\nu')$ as the cluster obtained
from $\K_n$ by adding the points in $K$ and not in $K_n$
with virtual multiplicities zero. Clearly,
$\T$ is consistent, $H_{\K_n}=H_{\T}$ and $\LL_{\T}=\sum_{p\in \K_+^Q}\alpha_p \LL_p$. 
It follows from \ref{Cart} that if $C$ is a curve going sharply through $\T$, then the strict transform $\widetilde{C}$ on $X$ is Cartier and $|\widetilde{C},L_p|_X=\alpha_p$ if $p\in
\K_+^Q$ and $0$ otherwise.
The blowing-up of $K'$ factors through $S_K$, so there is a birational morphism $g:S_{K'}\rightarrow
X$. Denote $T'_Q=\{p\in K'\mid g_*(E_p')=Q\}$. 
Write $\widetilde{C}^{K'}$ for the strict transform of $C$ on
$S_{K'}$. Then, by (\ref{for:int/excep}), we have that if
$p\notin T'_Q$, then
$|\widetilde{C}^{K'}\cdot
E_p^{K'}|_{S_{K'}}=\rho'_p=0$, so the direct image of $\widetilde{C}^{K'}$ by
$g$ intersects the exceptional locus of $X$ only at
$Q$. From this, it follows that $[\widetilde{C},L_p]_Q=|\widetilde{C},L_p|_X$ and this completes the proof.
\end{proofof}

An easy procedure for computing Cartier divisors with the
prescribed intersection multiplicities with the exceptional
components at $Q$ is derived from the proof of \ref{simpint}. 

\vspace{2mm}
\noindent\textbf{Procedure.} Keep the notation used there;  take $\K_0=BP(J)$. 

\noindent\textsc{Step 1}. Define $\K_1=(K_1,\nu^1)$ by adding to $\K_0$ a simple and free point $w$  in the
first neighbourhood of some point in $T_Q$ and unloading
multiplicities. 

\noindent\textsc{Step k-th}. While $\rho^{k-1}_p>0$ for some
$p\in \K_+^{Q}$, choose any such point $p$ and define $\K_k=(K_k,\nu^k)$ to be the cluster obtained by adding the minimal point not in $K_k$,
proximate to $p$ and infinitely near to $s_p$, and unload
multiplicities if the resulting cluster is not consistent. 

After finitely many steps, we reach the cluster $\K_n$. 
Define the cluster $\T$ as above. The strict transform on $X$ of a generic curve $C$  going through $\T$ intersects the exceptional divisor of $X$ only at $Q$ and with the prescribed interserction multiplicities. 

\vspace{2mm}
We close with an example.

\begin{example}
\label{ex:3} Let $I\in \I$ be an ideal with base points as in the Enriques diagram of figure \ref{fig:1} (Enriques diagrams are explained in \cite{EChisini} Book IV, Chapter 1 and also in  \cite{CasasBook} \S 3.9.). The dicritical points of $\K=BP(I)$ are $p_1,p_4,p_8$ and $p_{10}$ and so, the surface $X=Bl_I(S)$ has exceptional components $L_{p_1},L_{p_4},L_{p_8}$ and $L_{p_{10}}$. There is only one singularity on $X$, say $Q$. Take $\alpha_{p_1}=4,
\alpha_{p_4}=2, \alpha_{p_8}=4, \alpha_{p_{10}}=1$. Keeping
the notation as above, write
$J=I^4_{p_1}I^2_{p_4}I^4_{p_8}I_{p_{10}}$. The Enriques
diagrams of the clusters $\K$ and $\T$ are
shown in figure \ref{fig:1}. For any curve $C$ going
sharply through $\T$, $\widetilde{C}$ is a
Cartier divisor on $X$ locally irreducible near $Q$ and
$[\widetilde{C},L_{p_1}]_Q=4$,
$[\widetilde{C},L_{p_4}]_Q=2$,
$[\widetilde{C},L_{p_8}]_Q=4$ and
$[\widetilde{C},L_{p_{10}}]_Q=1$.
\begin{figure}
\begin{center}
 \psfrag{x}{$p_1$}\psfrag{b}{$p_2$}\psfrag{c}{$p_3$}\psfrag{d}{$p_4$}\psfrag{e}{$p_5$}
 \psfrag{f}{$p_6$}\psfrag{g}{$p_7$}\psfrag{h}{$p_8$}\psfrag{i}{$p_9$}\psfrag{j}{$p_{10}$}
\psfrag{q}{$O$}\psfrag{k}{$p_{r-1}$}\psfrag{l}{$p_r$}\psfrag{m}{$q_1$}\psfrag{n}{$q_2$}\psfrag{o}{$q_{r-1}$}\psfrag{p}{$q_r$}
\includegraphics[scale=0.5]{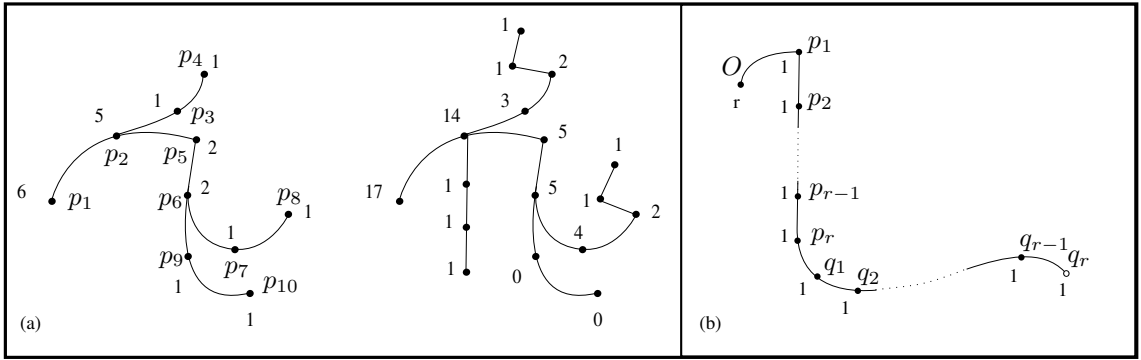}
\end{center}
\caption{\label{fig:1} (a) the Enriques diagrams of $\K$ and $\T$ in Example \ref{ex:3}; (b) the Enriques diagram $\mathbf{D}_r$, $r\geq 1$.}

\end{figure}
\end{example}

\subsection{The tangent space to the exceptional divisor}
Geometrical properties of the germ of the exceptional divisor of a birational projection of a normal singularity into a non-singular surface are derived here.

\begin{thm}
\label{tangent_space} Let $I\in \I$ and write $X=Bl_I(S)$ for the surface obtained by blowing-up $S$ along~$I$. Let $Q\in X$ be a singularity. 
The dimension of the tangent space to the reduced exceptional divisor $L$ at $Q$ is maximal, i.e. $\dim (\m_{L,Q}/\m_{L,Q}^2)=\sharp \K_+^Q$, where $\m_{L,Q}$ is the maximal ideal in the local ring $\OO_{L,Q}$. 
In particular, $\sharp \K_+^Q\leq \dim (\m_{X,Q}/\m_{X,Q}^2)$.
\end{thm}

\begin{rem}The reader may note that the no tangency of the exceptional components meeting at some singularity $Q\in X$ can easily be deduced from \ref{excessinchains}. The result stated here is stronger. 
\end{rem} 

To prove \ref{tangent_space} we need an easy lemma.

\begin{lem}
\label{lines} Let $\xi_1,\ldots,\xi_m$ be smooth curves in
$\mathbb{C}^N$ ($m\leq N$) going thought the point
$O=(0,0,\ldots,0)$ and for each $i$, write $l_i$ for the
tangent to $\xi_i$ at $O$. Assume there is a hypersurface
$H$ of $\mathbb{C}^N$ such that $[H,\xi_1]_O=1$ and
$[H,{\xi_i}]_O\geq 2$ for $i\in \{2,\ldots,m\}$. Then,
$l_1$ does not belong to the linear space generated by
$l_2,\ldots,l_m$ at $O$.
\end{lem}

\begin{proof}
If $H$ is a hypersurface of $\mathbb{C}^N$ such that
$[H,\xi_1]_O=1$, $H$ is necessarily smooth at $O$. If
moreover $[H,{\xi_i}]_O\geq 2$ for $i\in \{2,\ldots,m\}$,
the tangent space to $H$ at $O$ contains $l_2,\ldots,l_m$.
Therefore, if $l_1$ is contained in the linear space
generated by $l_2,\ldots, l_m$, it is also contained in the
tangent space to $H$ at $O$ and $[H,\xi_1]_O\geq 2$ against
the assumption.
\end{proof}

\begin{proofof}{\ref{tangent_space}}
Write $\K_+^Q=\{p_1,\ldots,p_m\}$. By virtue of
\ref{simpint}, there are Cartier divisors $\{C_i\}_{i\in
\{1,\ldots,m\}}$ on $(X,Q)$ such that
\[[C_i,L_{p_j}]_Q=\left \{
\begin{array}
{cc} 1 & \mbox{ if }i=j
\\ 2 & \mbox{ if }i\neq j. \end{array} \right.\]
Consider an embedding of $(X,Q)$ in a nonsingular ambiental
variety $\mathbb{A}_{\mathbb{C}}^{N}$. Since $C_i$ is
Cartier, it is the intersection of some hypersurface
$H_i:g_i=0$ in $\mathbb{A}_{\mathbb{C}}^{N}$ with $X$, and
$[H_i,L_{p_j}]_O=1$ if $i=j$ and 2 if$i\neq j$.
By \ref{lines} the tangent lines $\{l_{p_j}\}_{j=1,\ldots,m}$ to
the exceptional components $\{L_{p_j}\}_{j=1,\ldots,m}$
span a linear space of dimension $m=\sharp \K_+^Q$
and the first claim is proved. 
The second statement follows immediately. 
\end{proofof}

\emph{Minimal} singularities of a variety $V$ over $\mathbb{C}$ were introduced by Koll{\'a}r (Definition 3.4.1 of \cite{Kollar}):
$P\in V$ has a minimal singularity if ${\mathcal O}_{V,P}$ is reduced,
Cohen-Macaulay, the tangent cone of $V$ at $P$ is reduced and
\begin{eqnarray*}
\label{defn_min} \mult_P(V)=\emdim_P(V)-\dim_P(V)+1
\end{eqnarray*}
where $\emdim_P(V)$ is the embedding dimension
of $(V,P)$. 

\begin{cor}\label{excep_minimal}
The reduced exceptional divisor has only minimal singularities. 
\end{cor}

\begin{proofof}{\ref{excep_minimal}}
Let $Q\in X$ be a singular point and $L_p$ any exceptional component on $X$ going through $Q$. Then, from \ref{simpint}, we infer that $L_p$ is smooth at $Q$. 
Now, by virtue of Lemma 3.4.3 of \cite{Kollar} a curve singularity is minimal if and only if it has smooth branches with independent tangencies. The claim follows directly from \ref{tangent_space}.
\end{proofof}

The following corollary is immediate from \ref{tangent_space}.

\begin{cor}\label{nBP2}
If $(X,Q)$ is a normal surface singularity, there exists no
birational projection of it into a non-singular surface
contracting $n$ branches if their tangent directions at $Q$
are contained in a linear space of dimension $m<n$.
\end{cor}

\begin{rem}\label{bound:mult}
From \ref{tangent_space}, it turns out that  the number of exceptional components meeting at some singularity is upped bounded by the embedding dimension of the singularity. From Theorem 4 of \cite{Art66}, we infer that 
\begin{eqnarray}\label{1st_bound}
\mult_Q(X)\geq \sharp \K_+^Q-1.
\end{eqnarray}
In forthcoming \ref{nBP1}, we will characterize when the equality holds. However, this bound is weak in general as shown in the following example. In forthcoming \ref{bound_nQ}, we will obtain a sharper one.
\end{rem}

\begin{example}
\emph{Primitive} singularities are those sandwiched singularities that can be obtained by blowing-up a simple complete ideal (see \cite{Spivak} Definition 3.1). Here we construct primitive singularities $(X,Q)$ with multiplicity arbitrary high such that the exceptional locus of any birational projections at $Q$ is irreducible.
For any positive integer $r\geq 1$, we consider the Enriques diagram $\mathbf{D}_r$ defined by taking the origin $O$, $r$ consecutive vertices $p_1,\ldots, p_r$ all of them proximate to $O$ and $r$ consecutive free vertices $q_1,\ldots,q_r$ infinitely near to $p_r$ (see (b) in figure \ref{fig:1}). Take any simple ideal $I_p$ such that $\mathbf{D}_r$ is the Enriques diagram of $\K=BP(I_p)$. Then,  the exceptional divisor of the surface $X$ obtained by blowing up $I_p$ is irreducible, so $\sharp \K_+^Q=1$. However, it can easily be seen that the multiplicity of $(X,Q)$ is $r+1$ (use for example, Theorem 4.7 of \cite{moi1}). 
\end{example}

\begin{rem}
Sandwiched singularities can be constructed by means of the so-called \emph{decorated curves} introduced in \cite{JS}. 
The bound stated in (\ref{1st_bound}) follows from Remark 3.6 of \cite{deJong} for sandwiched singularities $X(C,l)$ provided that $l\gg0$ (see the notation used in \cite{deJong}). However, our argument works for sandwiched singularities with no restrictions.
\end{rem}

\vspace{2mm}
\section{A bound for the number of exceptional components} 

Here, we establish a formula for the number of branches of a hypersurface section through a sandwiched singularity. 
From it, we will obtain an upper bound for the number of exceptional components meeting at the same singularity. 

First of all, we need a couple of easy lemmas and to introduce some notation. 
\begin{lem}\label{for:multQ}
$\mult_Q(X)=1+\sharp \B_Q$.
\end{lem}

\begin{proof}
By theorem 4.7 of \cite{moi1}, $\mult_Q(X)=\K_Q^2-\K^2$, where $\T^2$ means the self-intersection of the cluster $\T=(T,\tau)$, $\T^2=\sum_{q\in T}\tau_q^2$. The assertion follows by direct computation using that $I_Q$ has codimension one in $I$ and the Proposition 4.7.1 of \cite{CasasBook}.
\end{proof}

Write $\K_+^{Q}=A_1\sqcup A_2$ where
$A_1=\{p\in \K_+^{Q} \mid p\geq o_{Q}\}$ and
$A_2=\{p\in \K_+^{Q} \mid p \ngeq o_{Q}\}$.
\begin{lem}
\label{lemaA} We have that \rm{(i)} $A_1=\K_+^Q\cap \B_{Q}$; \rm{(ii)} $A_2\subset \{p\in K \mid o_Q\rightarrow p\}$.
\end{lem}

\begin{proof}
(i) The inclusion $\supset$ follows from (c) of \ref{coef_fund}. Now, let $p\in \K_+^Q$. Then $p$ must be is adjacent to some point in $T_Q$. If $p\geq o_Q$, $p$ is $m_K$-proximate to some point in $T_Q$. (c) of \ref{coef_fund} gives the claim. 

(ii) Let $p\in A_2$. By \ref{difexcess}, we have that $\varepsilon_p-\sum_{q\rightarrow p}\varepsilon_q=-1$. Since $p\ngeq o_Q$, $\varepsilon_p=0$ and $\sum_{q\rightarrow p}\varepsilon_q=1$. It follows that $o_Q$ is proximate to $p$. 
\end{proof}

A point $q$ is said to be \emph{$T_Q$-free} (respectively, \emph{$T_Q$-satellite}) if it is
proximate to one point $p$ in $T_Q$ (resp. two points $p_1,p_2\in T_Q$). Notice that (c) of \ref{coef_fund} says that $\B_Q\subset \{q \rightarrow T_Q\}$. Thus, we can split $\B_Q$ in $\B_Q=\B_Q^1\cup \B_Q^2$, where $\B_Q^1=\{q\in \B_Q \mid q\mbox{ is $T_Q$-free}\}$ and $\B_Q^2=\{q\in \B_Q \mid q\mbox{ is $T_Q$-satellite}\}$. The reader may note that if $q$ is $T_Q$-satellite, then $q\in T_Q$, so
\begin{eqnarray}\label{auxB2}
\B_Q^2\cap T_Q=\B_Q^2.
\end{eqnarray} 

\begin{prop}
\label{branches} The number of branches of a generic  hypersurface section of
$(X,Q)$ is
\[\br=\mult_Q(X)-\sharp \B_Q^1\cap T_Q.\]Moreover, $\br\geq \K_+^Q-1$ and the equality holds if and only if the following three conditions are satisfied:
\begin{itemize}
 \item[(i)] $\B_Q^2=\emptyset$;
\item[(ii)] every point in $\B_Q\setminus T_Q$ is $m_K$-proximate to $T_Q$; 
\item[(iii)] $o_Q$ is $m_K$-satellite.
\end{itemize}
\end{prop}

\begin{proof}
First of all, since $I_Q\OO_X={\mathcal M}_QI\OO_X$ (see \S 2.2), the strict transform of generic curves going through $\K_Q$ are generic hypersurface sections of $(X,Q)$. Thus, $\br$ 
equals the intersection number of the strict transform
of a curve $C$ going sharply through $\K_Q$ and the reduced exceptional divisor in the minimal resolution of $(X,Q)$, i.e.
\begin{eqnarray} \label{numberbranches}
\br=\sum_{p\in T_Q}|\widetilde{C}^K{\cdot}
E_p|_{S_K} .\end{eqnarray}
From (\ref{for:int/excep}) and (a) of
\ref{difexcess}, this equals to \[\sum_{p\in T_Q}\rho'_p=\sum_{p\in T_Q}
\varepsilon_p+\sum_{p\in T_Q}\sharp\{q\in \B_Q \mid q
\rightarrow p\}.\]The definition of the
$\varepsilon$'s gives that $\sum_{p\in T_Q}\varepsilon_p
= 1-\sharp \B_Q\cap T_Q$. On the other hand, 
\[\sum_{p\in T_Q}\sharp\{q\in \B_Q \mid q\rightarrow p\}=\sharp \B_Q^1+2\,\sharp \B_Q^2= \sharp\B_Q+\sharp\B_Q^2.\]It follows from this and (\ref{auxB2}) that 
$\br=1+\sharp \B_Q-\sharp \B_Q^1\cap T_Q$.  \ref{for:multQ} completes the proof of the first assertion.
Notice that $A_1\subset \B_Q\setminus \B_Q\cap T_Q$. Thus, by virtue of \ref{lemaA}, we have that 
$\br=\sharp A_1+ \sharp A_2 \geq \sharp \B_Q -\sharp \B_Q\cap T_Q +2$. 
It follows that $\br\geq \sharp \K_+^Q-1+\sharp \B_Q^2\geq \K_+^Q-1$ and the equality holds if and only if the above three conditions hold. 
\end{proof}

Part of the interest of the preceding result lies on the fact that it allows to know if $\K_+^Q=\br+1$ by checking conditions (i)-(iii) directly in the Enriques diagram of the base points of $I$. 
\begin{rem}
A natural question is if for any sandwiched singularity $(X,Q)$ there exists a birational projection to a non-singular surface such that $\br=\K_+^Q-1$. The anwer is no in general as shown in the following example (cf. Theorem \ref{nBP1}).
\end{rem}

\begin{example}\label{ex:4}
Take $I\in \I$ with base points as shown in Figure \ref{fig:4}. By blowing-up $I$ we obtain a surface $X$ with just one sandwiched singularity, say $Q$. The exceptional divisor of $X$ is irreducible, so $\K_+^Q=1$. It can be seen by using Theorem 3.2 of \cite{AlbFer} that the proximities between the vertices of $\G_Q$ are in this case determined. That is, assume that $J\in \I$ is a complete ideal such that there is a singularity $Q'=X'=Bl_{J}(S)$ such that $\OO_{X',Q'}$ is analytically isomorphic to $\OO_{X,Q}$; then there is a bijection $\psi:T_{Q'}\rightarrow T_Q$ such that $p\rightarrow q$ if and only if $\psi(p)\rightarrow \psi(q)$ (with the language of \cite{AlbFer}, we say that there is only one contraction for $(X,Q)$). In particular, $p\in \B_{Q'}^2$ if and only if $\psi(p)\in \B_Q^2$.  Notice that $p_4,p_5\in \B_Q^2$. Thus, there is no ideal for $(X,Q)$ satisfying the assumpion (i) of \ref{branches}.
\begin{figure}
\begin{center}
 \psfrag{a}{$p_1$}\psfrag{b}{$p_2$}\psfrag{c}{$p_3$}\psfrag{d}{$p_4$}\psfrag{e}{$p_5$}
 \psfrag{f}{$p_6$}\psfrag{g}{$p_7$}
\includegraphics[scale=0.5]{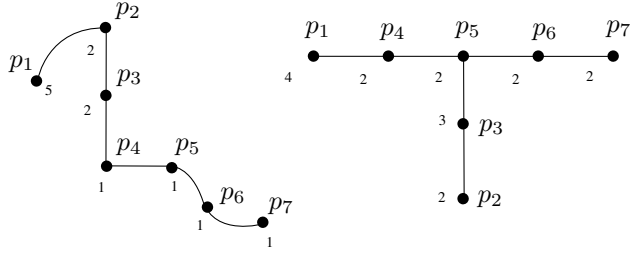}
\end{center}
\caption{\label{fig:4} The Enriques diagrams of the base points of the ideal $I$ and the dual graph of $(X,Q)$ in Example \ref{ex:4}.}

\end{figure}
\end{example}

\begin{cor}\label{bound_nQ}
If $(X,Q)$ is a normal surface singularity, there exists no
birational projection of it into a non-singular surface
contracting more that $\br+1$ smooth branches through $Q$. 
\end{cor}

\vspace{2mm}

To close this section, we characterise when
the bound in (\ref{1st_bound}) is attained. 
The following lemma provides some technical characterisations for minimal singularities. 

\begin{lem}
\label{minimal_1} The singularity $(X,Q)$ is minimal if and only if one of the following equivalent conditions is satisfied:
\begin{itemize}
\item[(i)] the fundamental cycle $Z_Q$ is reduced; \item[(ii)] $\br=\mult_Q(X)$; \item[(iii)]
$\B^1_Q\cap T_Q=\emptyset$.
\end{itemize}
\end{lem}

\begin{proof} 
By virtue of 3.4.10 of  \cite{Kollar}, a normal surface singularity $(X,Q)$ is  minimal if and only if it is sandwiched and the fundamental cycle $Z_Q$ is reduced (cf. \S2 of \cite{Spivak}, Lemma 5.8 of \cite{Bondil-Le}).
Thus (i) is equivalent to the minimality of $(X,Q)$. 

\vspace{2mm} \noindent(i)$\Leftrightarrow$(ii) Let $C\subset S$
be a curve going sharply through $\K_Q$ so that
$\widetilde{C}$ is a transverse hypersurface section of $(X,Q)$.  The projection formula applied to $f$ gives that $\mult_Q(X)=\sum_{p\in
T_Q}z_p|\widetilde{C}^K{\cdot}E_p|_{S_K}$. If $Z_Q$ is reduced, (ii) follows from (\ref{numberbranches}) above.

\vspace{2mm} \noindent(ii)$\Leftrightarrow$(iii) follows directly from the formula of \ref{branches}.
\end{proof}

\begin{prop} \label{excep=bound}
Let $Q\in X$ be a singular point. Then, we have
$\sharp \K_+^Q\leq \emdim_Q(X)$, and the equality holds if
and only if the three following conditions are satisfied:
 \begin{itemize}
 \item[(1)] $\B_Q\cap T_Q=\emptyset$;
 \item[(2)] every point in $\B_Q$ is $m_K$-proximate to $T_Q$;
 \item[(3)] $o_Q$ is $m_K$-satellite
 \end{itemize}
In particular, a necessary condition for the equality to hold is that $(X,Q)$
is minimal.
\end{prop}

\begin{rem}
The reader may note that the conditions (1)-(3) imply immediately the conditions (i)-(iii) of \ref{branches}.
\end{rem}

\begin{proofof}{\ref{excep=bound}}
The inequality $\sharp \K_+^Q\leq \emdim_Q(X)$ has been proved in \ref{tangent_space}. 
By \ref{lemaA}, we have that
$A_1=\K_+^Q\cap \B_Q$. It follows that $A_1 \subset
(\B_Q\setminus T_Q)\subset \B_Q$, and the equality holds if and only if $\B_Q\cap T_Q=\emptyset$ and every $q\in\B_Q$ is $m_K$-proximate to $T_Q$. 
Similarly, we know that $A_2\subset \{p\in K \mid
o_{w}\rightarrow p\}$. Therefore,
$\sharp \K_+^Q = \sharp A_1 +\sharp A_2 \leq \sharp \B_Q+ 2$
and the equality is satisfied if and only if 
(1) $\B_Q\cap T_Q=\emptyset$, 
(2) every  $q\in \B_Q$ is $m_K$-proximate to $T_Q$ and 
(3) $o_Q$ is $m_K$-satellite.
Finally, note that by \ref{minimal_1}, condition (1) implies immediately that $(X,Q)$ is minimal.
\end{proofof}

From the above result we obtain that minimal singularities are just those normal singularities that admit a birational projections contracting $\mult_Q(X)+1$ smooth branched through it. 

\begin{thm}\label{nBP1}
Let $(X,Q)$ be a normal surface singularity. Then, $(X,Q)$ is minimal if and only if there exists a birational projection contracting exactly $\mult_Q(X)+1$ (smooth) branches through~$Q$.
\end{thm}

\begin{proof}
The ``if'' part follows from \ref{excep=bound}. The ``only if'' part is a consequence of the fact that for every minimal singularity, there exists a a complete ideal $J$ such that the cluster of its base points $BP(J)$ satisfies the conditions (1)-(3) of \ref{excep=bound} and $(X,Q)$ lies on the surface obtained by blowing-up $J$. 
Indeed, let $S'$ be the surface obtained by blowing-up the point $O\in S$ and a point $u$ in the first neighbourhood of $O$; let $R'$ be the local ring $\OO_{S',p}$ where $p$ is the (satellite) point in the intersection of the two exceptional components of $S'$. 
Now, an ideal $J'\in \mathbf{I}_{R'}$ can be considered so that it has only free base points, it satisfies the conditions of Corollary 1.14 of \cite{Spivak} relative to the singularity $(X,Q)$ and no base point of $J'$ in the first neighbourhood of $p$ is proximate to $O$ or $u$. 
Write $\K'=(K',\nu')$ for the cluster of base points of $J'$. By blowing-up $J'$ we obtain a surface $X'=Bl_{J'}(S')$ with just one singular point $Q'$ so that $\OO_{X,Q}$ is analytically isomorphic to $\OO_{X',Q'}$. 
Now, for each $p\in T_Q$, write $\v(p)$ for the weight of $p$ in the resolution graph $\G_Q$ and $\deg(p)$ for the number of vertices in $\G_Q$ adjacent to $p$. Then, it is immediate to see that $p=o_Q$ so $\sharp \{u\in K' \mid u\rightarrow p\}=\deg(o_Q)$, while for $q\in T_Q$, $q\neq o_Q$, we have $\sharp \{u\rightarrow q\} =\deg(q)-1$. It follows that the number of dicritical vertices of $\K'$ equals 
\begin{eqnarray*}
 \sum_{q\in T_Q} (\omega(q)-\deg(q))-1 = \mult_Q(X)-1
\end{eqnarray*}
the last equality because $Z_Q$ is reduced.
On the other hand, the surface $X'$ dominates $S'$ and so, it also dominates $S$. 
Then, it is enough to take the projection of $X'$ to $S$ in order to obtain a birational map contracting $\mult_Q(X)+1$ smooth branches. 
\end{proof}

\begin{rem}
Notice that the conditions (1)-(3) of \ref{excep=bound} can be checked in the Enriques diagram of the cluster $\K$. Thus the ``only if'' part of the above theorem can also be proved by constructing by hand an Enriques diagram for $(X,Q)$ satisfying (1)-(3) and making use of Proposition 2.1 of \cite{AlbFer} (see \cite{AlbFer} for details).
\end{rem}

\bibliographystyle{amsplain}

\end{document}